\magnification\magstep1
\baselineskip = 18pt
\def\n{\noindent}
\def \reel{ {\rm I}\!{\rm R}}

\def \comp{ \;{}^{ {}_\vert }\!\!\!{\rm C}   }

\def \nat{ { {\rm I}\!{\rm N}} }

\def\pf{\medskip{\noindent{\bf Proof. }}}
\def \rat{ {\rm Q}\kern-.65em {}^{{}_/ }}

 \baselineskip = 18pt
 \def\n{\noindent}
 \overfullrule = 0pt
 \def\qed{{\hfill{\vrule height7pt width7pt
depth0pt}\par\bigskip}} 
\def \reel{ {\rm I}\!{\rm R}}

\centerline{\bf Regular operators}
\centerline{\bf between non-commutative $L_p$-spaces.}\bigskip
\centerline{by}\bigskip
\centerline{Gilles Pisier\footnote*{Supported in part by the NSF}}
\centerline{Texas A\&M University}
\centerline{and}
\centerline{Universit\'e Paris VI}\vskip12pt

{\bf Abstract}
We introduce the notion of a regular mapping on a 
non-commutative
$L_p$-space associated to a hyperfinite
 von Neumann algebra for $1\le p\le \infty$. This is 
a non-commutative generalization of the 
notion of regular or order bounded map on a Banach lattice.
This extension is based on our
recent paper [P3], where we introduce and study 
a non-commutative version of
vector valued $L_p$-spaces. In the extreme cases $p=1$ and $p=\infty$, our regular operators
reduce to the completely bounded ones and the regular norm coincides
with the $cb$-norm. We prove that a mapping is regular
iff it is a linear combination of 
bounded, completely positive mappings. We  prove an extension
theorem for regular mappings defined on a subspace
of a non-commutative
$L_p$-space. Finally, let $R_p$ be the space of all regular
mappings on a given non-commutative
$L_p$-space  equipped with 
the regular norm. We prove the isometric identity
$R_p=(R_\infty,R_1)^\theta$ where $\theta=1/p$ and where
$(\ .\ ,\ .\ )^\theta$ is the dual
 variant of Calder\'on's complex interpolation method.

\vskip24pt
\item{} {\bf Plan}
\item{} Introduction.
\item{\S 1.} Background on operator spaces.
\item{\S 2.}Completely positive maps and regular maps.
\item{\S 3.}Main result and applications.\vfill\eject

\n {\bf \S 0. Introduction.}

Let $L,\Lambda$ be Banach lattices. An operator $u\colon \ L\to \Lambda$ is
called regular (or sometimes ``order bounded'') if there is a constant $C$
such that for all finite sequences $x_1,\ldots, x_n$ in $L$ we have
$$\|\sup_{i\le n}|u(x_i)|\ \|_\Lambda\le C\|\sup_{i\le n}|x_i|\|_L.\leqno
(0.1)$$
We denote by $\|u\|_r$ the smallest constant $C$ such that this holds. More
generally, this definition makes sense for an operator $u\colon \ S\to
\Lambda$ defined only on a subspace $S\subset L$ by considering (0.1)
restricted to $n$-tuples of elements of $S$.

We will denote by $B_r(L,\Lambda)$ (resp. $B_r(S,\Lambda))$ the space of all
such operators $u\colon \ L\to \Lambda$ (resp. $u\colon\  S\to \Lambda$)
and we equip it with the norm $\|\ \|_r$. Equivalently, let us denote
$S[\ell^n_\infty]$ the space $S\otimes \ell^n_\infty$ equipped with the
norm $\left\|\sum\limits^n_{i=1} x_i\otimes e_i\right\| =
\|\sup\limits_{i\le n} |x_i|\ \|_L$. Then an operator $u\colon \ S\to
\Lambda$ is regular iff there is a constant $C$ such that for all $n$
$$\|u\otimes I_{\ell^n_\infty}\|_{S[\ell^n_\infty]\to
\Lambda[\ell^n_\infty]} \le C.\leqno (0.2)$$
Moreover $\|u\|_r$ is equal to the smallest constant $C$ for which this
holds. It is well known (see [MN, p. 24]) that if $\Lambda$ is 
Dedekind complete (this is
the terminology of [MN], the
same notion is sometimes called order complete) then
$u\in B_r(L,\Lambda)$ iff $u$ is a linear combination of bounded positive
maps from $L$ into $\Lambda$. This is classical. More recently, we observed
([P2]) that the space $B_r(L_p,L_p)$ coincides with the complex
interpolation
space $(B(L_\infty,L_\infty), B(L_1,L_1))^\theta$ when $\theta = 1/p$. We
also proved an extension theorem for $L_p$-valued regular maps defined on a
subspace of $L_p$. In this paper, we wish to extend the notion of
regularity and to present a non-commutative version of these results of
[P2].
 In
particular, we will consider the Schatten $p$-class $S_p$, a subspace
$S\subset S_p$ and an operator $u\colon\ S\to S_p$. It turns out that
(0.1) and (0.2) can be formulated in this setting (although of course $S_p$
is not a Banach lattice in the usual  sense).

 In the non-commutative case, the space $S[\ell^n_\infty]$ is replaced by
$S[M_n]$ where $M_n$ is the Banach space of all complex $n\times n$ matrices,
equipped with the norm of $B(\ell^n_2)$. This extension is based on our
recent paper [P3] (see also the announcement [P5]), where we introduce and study a non-commutative version of
vector valued $L_p$-spaces. In that theory, the space of values is an
operator space, i.e. a closed subspace $E\subset B(H)$ for some Hilbert
space $H$. In [P3] we defined the $E$-valued version of the Schatten
$p$-class and we denoted it by $S_p[E]$. If $S\subset S_p$ is a closed
subspace, we  denote by $S[E]$ the closure of $S\otimes E$ in $S_p[E]$. In
particular this makes sense if $E=M_n$ and it allows us to extend the
definition of regular operators.

 In particular, we  will prove below

\proclaim Theorem 0.1. Let $1\le p <\infty$. Let $u\colon \ S_p\to S_p$ be
a linear map. The following are equivalent.\medskip
\item{(i)} There is a constant $C$ such that for all $n$, $u\otimes
I_{M_n}$
is bounded on $S_p[M_n]$ and\break $\|u\otimes I_{M_n}\|_{S_p[M_n]\to
S_p[M_n]} \le C$.
\item{(ii)} The map $u$ is a linear combination of bounded, completely
positive maps on $S_p$.
\item{(iii)} There are completely positive, bounded linear maps $u_j\colon
\ S_p\to S_p$ $(j=1,2,3,4)$ such that $u = u_1-u_2+i(u_3-u_4)$.\medskip

\n Here of course a linear map $v\colon \ S_p\to S_p$ is called completely
positive if for each $n\ge 1$ the map $I_{M_n}\otimes v\colon \
S_p(\ell^n_2\otimes \ell_2)\to S_p(\ell^n_2\otimes\ell_2)$ is positive
(i.e. positivity preserving). This extends the usual complete  positivity
in the case $p=\infty$.

We will abbreviate completely positive by $c.p.$.
\n  Note that in [P1] we
proved the following result (see [P1, section~8]).

\proclaim Theorem 0.2. Let $0<\theta<1$ and $p=1/\theta$. Then the space
$$(cb(S_\infty, B(\ell_2)), \quad cb(S_1,S_1))^\theta$$
coincides with the space of all operators $u\colon \ S_p\to S_p$ which are
linear combinations of bounded $c.p.$ maps on $S_p$ (or equivalently which
satisfy (iii) in Theorem~0.1).

In section~1 below, we summarize the results from [P1, P3] that we will
need. In section~2 we discuss completely positive maps and regular maps on
$S_p$. In section~3, we prove our main results including Theorem~0.1.
Finally we will prove an extension theorem:\ any regular map $u\colon \
S\to S_p$ defined on a subspace of $S_p$ extends to a regular map on $S_p$
with the same regular norm. This generalizes to $S_p$ the well known
extension property of $c.b.$ maps into $B(\ell_2)$.

We have restricted ourselves for simplicity to the simplest non-commutative
$L_p$-space, i.e. $S_p$, but it is clearly possible to extend our results
to the setting of  $L_p$-spaces associated to a  hyperfinite
 (or equivalently injective, by Connes's well
known results [Co])
von~Neumann algebra $M$ equipped with a 
semi-finite trace as in [P3]. See  [N, H2, Ko, Te1, Te2]
and the references in [N], 
for the classical ({\it i.e.}
 scalar valued) theory of non-commutative $L_p$-spaces.  The resulting statements
can be viewed as an extension  to a non-commutative $L_p$-space
	 of 
Wittstock's well known extension/decomposition theorem
for $cb$ maps into an injective 
von Neumann algebra. Indeed, viewing  a von Neumann algebra
as a 
non-commutative $L_\infty$-space, Wittstock's result corresponds
to the case $p=\infty$. We state this extension 
without proof at the end of the paper.
We leave
the details to the reader.

\vfill\eject
\n {\bf \S 1. Background on operator spaces.}

In this section, we recall a number of notions and results that we will
use.

By an operator space, we mean a closed subspace $E\subset B(H)$ of the space
$B(H)$ of all bounded operators on some Hilbert space $H$. We will use
the theory of operator spaces, as developed recently in a series of papers
[BP, ER2]. In this theory, an operator space structure is a sequence of
norms on the spaces $M_n(E)$ of all $n\times n$ matrices with entries in
$E$. If $E\subset B(H)$ is a subspace of $B(H)$, its natural operator space
structure is by definition the one corresponding to the norms induced on
$M_n(E)$ by $B(\ell^n_2(H))$, more precisely to the norms $\| \ \|_n$
defined on $M_n(E)$ by
$$\|(a_{ij})\|_n = \|(a_{ij})\|_{B(\ell^n_2(H))}.\leqno (1.1)$$
Given two embeddings $E\subset B(H)$, $E\subset B(K)$ of the same space
into $B(H)$ and $B(K)$ we will say that they define the same operator space
structure if the associated norms (1.1) are the same for both embeddings.

By a fundamental result of Ruan [Ru] the sequences of norms $\|\ \|_n$
which arise in this way can be characterized by very simple axioms (see
[Ru]). We will refer to the sequence (1.1) of norms $\|\ \|_n$ on $M_n(E)$
as ``the operator space structure'' on $E$. Note that the Banach space
structure of $E$ is included via the case $n=1$.

Thus, to give ourselves an operator space structure on a vector space $E$
is the same as to give ourselves an equivalence class of embeddings $E\to
B(H)$ (two embeddings are equivalent if they give the same structure via
(1.1)). Let $E,F$ be operator spaces, a map $u\colon \ E\to F$ is called
completely bounded (in short $c.b.$) if $I_{M_n}\otimes u$ is bounded from
$M_n(E)$ into $M_n(F)$ and $\sup\limits_{n\ge 1} \|I_{M_n}\otimes
u\|<\infty$. We define $\|u\|_{cb} = \sup\limits_{n\ge 1} \|I_{M_n}\otimes
u\|$. We denote by $cb(E,F)$ the space of all such maps equipped with this
norm. We refer to [Pa] for more information. We will say that $u\colon\
E\to F$ is completely isomorphic (resp. completely isometric) if $u$ is an
isomorphism and $u$ and $u^{-1}$ are $c.b.$ (resp. if $I_{M_n}\otimes u\colon \
M_n(E) \to M_n(F)$ is isometric for all $n\ge 1$).

In [BP, ER2], it was proved that $cb(E,F)$ can be equipped with an operator
space structure by giving to $M_n(cb(E,F))$ the norm of the space $cb(E,
M_n(F))$. In particular, this defines an operator space structure on the
dual $E^* = cb(E,\comp)$ so that we have isometrically
$$M_n(E^*) =
cb(E,M_n). \leqno (1.2)$$
This duality for operator spaces
 has many of the nice properties of the
usual Banach space  duality.
 See [BP,ER2,B2,BS,ER4] for more details.

  In [P1], we introduced complex interpolation
for operator spaces. Let $(E_0,E_1)$ be a compatible couple of Banach spaces
in the sense of interpolation theory. We refer to [BL,Ca] for the basic
notions and the definitions of the complex interpolation methods
$(E_0,E_1)_\theta$ and $(E_0,E_1)^\theta$.

\n Now assume $E_0,E_1$ each equipped with an operator space structure (in the
form of norms on $M_n(E_0)$ and $M_n(E_1)$ for all $n$). Let $E_\theta =
(E_0,E_1)_\theta$ and $E^\theta = (E_0,E_1)^\theta$. Then, we can define
an operator space structure on $E_\theta$ (resp. $E^\theta$) by setting
$$M_n(E_\theta) = (M_n(E_0), M_n(E_1))_\theta \ (\hbox{resp. }
M_n(E^\theta) = (M_n(E_0), M_n(E_1))^\theta).\leqno (1.3)$$
In [P1] we observed that these norms verify Ruan's axioms and hence they
define an operator space
 structure on $E_\theta$ (resp. $E^\theta$).

  Let
$H,K$ be Hilbert spaces. Let $1\le p<\infty$. Let us denote by $S_p$ (resp.
$S_p(H)$, $S_p(H,K)$) the space of all 
operators $T\colon \ \ell_2\to \ell_2$
(resp. $T\colon \ H\to H$, resp. $T\colon \ H\to K$) such that
${\rm tr}|T|^p<\infty$. We equip it with the norms
$$\|T\|_p = ({\rm tr}|T|^p)^{1/p}.$$
For notational convenience, we denote by $S_\infty$ (resp. $S_\infty (H),
S_\infty(H,K)$) the space of all compact operators on $\ell_2$ (resp. on
$H$, resp. from $H$ into $K$), equipped with the operator norm.

Since $S_\infty$ and $B(H)$ are $C^*$-algebras they have an
obvious operator
space structure, which we will call the {\it  natural} operator space structure
(in short o.s.s.). More generally, the space
 $B(H,K)$ has a {\it  natural} o.s.s.
corresponding to the norm induced by $B(\ell^n_2(H), \ell^n_2(K))$ on
$M_n(B(H,K))$. Since $S_1 = S^*_\infty$, we can equip $S_1$ with the dual
operator space structure, as defined above in (1.2). Similarly for any
predual $X$ of a von~Neumann algebra, we obtain by (1.2) an o.s.s. on $X$
which we will call the {\it natural\/} o.s.s. on $X$.

 In particular for any	measure space $(\Omega,\mu)$, the spaces
$L_\infty(\mu)$ and $L_1(\mu)$ can now be viewed as operator spaces with
their {\it natural\/} o.s.s. Let $0<\theta<1$ and $p=1/\theta$. By
interpolation, using
$$L_p(\mu) = (L_\infty(\mu), L_1(\mu))_\theta$$
and
$$S_p = (S_\infty, S_1)_\theta$$
we define using (1.3) an o.s.s. on $L_p(\mu)$ and $S_p$, which we will
again call the {\it natural\/} one.

 We now turn to the vector valued case.
Let $E,F$ be operator spaces. We denote by $E\otimes_{\rm min} F$ the
completion of $E\otimes F$ for the minimal (= spatial) tensor product. If
$E\subset B(H)$, then $L_\infty(\mu; E) \subset L_\infty(\mu; B(H))$ so
that
we can equip $L_\infty(\mu; E)$ with the o.s.s. induced by $L_\infty(\mu;
B(H))$. Similarly, we define $S_\infty[E] = S_\infty \otimes_{\rm min} E$
and we equip it with the natural o.s.s. of the minimal (= spatial) tensor
product. In the finite dimensional case, we define
$S_\infty^n[E] = M_n \otimes_{\rm min} E$. This is the same as $M_n(E)$,
which we will sometimes also denote by $M_n[E].$
This settles the simplest case $p=\infty$. 

We now turn to the case
$p=1$.
We define
$$L_1(\mu;E) = L_1(\mu) \otimes^\wedge E$$
and
$$S_1[E] = S_1\otimes^\wedge E,$$
where the symbol $\otimes^\wedge$ refers to the operator space version of
the projective tensor norm introduced in [BP, ER2] and developed in [ER6].
This defines operator space structures in $L_1(\mu; E)$ and $S_1[E]$. By
interpolation, we define
$$L_p(\mu; E) = (L_\infty(\mu; E), L_1(\mu; E))_\theta \hbox{ and } S_p[E] =
(S_\infty[E], S_1[E])_\theta,\leqno (1.4)$$
and using (1.3) this gives us an operator space structure on $L_p(\mu; E)$
and $S_p[E]$. We will refer to all 
these structures in each case as 
{\it the natural  one}.

Clearly, a similar definition applies to $S_p(H)$ and $S_p(H,K)$. We will
denote by $S_p[H;E]$ and $S_p[H,K; E]$ the resulting operator spaces.
Moreover if $H = \ell^n_2$ we will denote simply by $S^n_p[E]$ the space
$S_p[\ell^n_2;E]$. Note that by definition $S_p[E]$ is included into the
space $M_\infty(E)$, and $S^n_p[E]$ (which is clearly linearly isomorphic
to $M_n(E)$) is completely isometric to the subspace of $S_p[E]$ formed by
the matrices with entries in $E$ which have zero coefficients outside the
``upper left corner'' of size $n\times n$.

We will also need an explicit formula for the norm in the spaces $S_p[H;E]$
or $S_p[H,K;E]$ as described in [P3].

For any $x$ in $S_p[H,K;E]$ (viewed as an element of the larger space
$S_\infty(H,K) \otimes_{\rm min} E$) we have
$$\|x\|_{S_p[H,K;E]} = \inf\{\|a\|_{S_{2p}(\ell_2 ,K)} \|y\|_{S_\infty[E]}
\|b\|_{S_{2p}(H,\ell_2)}\}\leqno (1.5)$$
where the infimum runs over all possible decompositions of $x$ of the form
$$x = (a\otimes I_E)\cdot (y)\cdot (b\otimes I_E),$$
where  the dot denotes the ``matrix'' product, and  we identify
$S_{2p}(\ell_2,K)$ and $S_{2p}(H,\ell_2)$ with the corresponding matrix
spaces and $S_\infty[E]$ is viewed as a subset of the set of all matrices
$(y_{ij})_{i,j\in {\nat}}$ with entries in $E$. In the particular case
when $H,K$ are finite dimensional, we have for all $x$ in $S^n_p[E]$
$$\|x\|_{S^n_p[E]} = \inf\{\|a\|_{S^n_{2p}} \|y\|_{M_n(E)}
\|b\|_{S^n_{2p}}\}\leqno (1.5)'$$
where the infimum runs over all possible decompositions
$$x = (a\otimes I_E) \cdot y \cdot (b\otimes I_E)$$
(matrix product).
Indeed, if $H,K$ are $n$-dimensional in (1.5) we can replace $a$ and $b$
respectively by $aP$ and $Qb$ where $P$ (resp. $Q$) is the orthogonal
projection onto the essential support of a (resp. onto the range of $b$);
if we also replace $y$ by $(P\otimes I_E)\cdot y \cdot (Q\otimes I_E)$ we
easily obtain $(1.5)'$.

We will use the following version of Fubini's theorem
(cf. [P3, Theorem 1.9]): 

\n  Let $1\le p \le \infty$. Let $H,K$ be arbitrary
Hilbert spaces and let $E$ be an operator space. We have
completely isometrically
$$S_p[H; S_p[K; E]] \simeq S_p[H\otimes_2 K; E] \simeq S_p[K; S_p[H; E]].\leqno(1.6)$$

Let $m\ge 1$ be an arbitrary integer. 
For any matrix $a$ in $M_m$, let us denote 
by $tr_m(a)$ its  trace in the usual sense.
Let $z$ be an
element of $S^m_1[M_m]$, of the form\break $z=(z_{ij})$ with
$z_{ij}\in M_m$. Then we have linear identifications
$S^m_1=M^*_m$ and \break $S^m_1[M_m] = S^m_1\otimes M_m$. These allow to
define the trace functional $z\to tr(z)$ associated to the
pairing $M^*_m \times M_m\to {\comp}$. It is easy to check
that $$tr(z) = \sum_{ij} \langle z_{ij}, e_{ij}\rangle  =
\sum_{ij} \langle e_i, z_{ij}(e_j)\rangle.$$ Moreover if
$z=a\otimes b$ with $a\in S^m_1$ and $b\in M_m$ then we
have $$tr(z) = tr_m(^tab)\leqno (1.7)$$ where $^ta$ is the
transposed of $a$.

We will use the following elementary identity valid for any $\alpha,\beta$
in $M_m$
$$tr((\alpha\otimes I)\cdot z\cdot (\beta\otimes I)) = tr((I\otimes
^t\alpha) \cdot z\cdot (I\otimes ^t\beta)).\leqno (1.8)$$
This is easily checked as follows. It clearly suffices to prove this for
$z=a\otimes b$. Then $(\alpha\otimes I)\cdot z\cdot (\beta\otimes I) =
(\alpha a\beta) \otimes b$ so that (1.8) reduces to
$$tr_m(^t(\alpha a\beta)b) = tr_m(^ta(^t\alpha b^t\beta))$$
which follows from the well known identity satisfied by the trace of the
product of two matrices.

Finally, we will use the inequality
$$|tr(z)| \le \|z\|_{S^m_1[M_m]}.\leqno (1.9)$$
This is an immediate consequence of the construction of the projective
tensor product $S^m_1\otimes^\wedge M_m$ as developed in [BP, ER2]. For the
convenience of the reader, we briefly sketch a direct deduction of (1.9)
from $(1.5)'$.

Assume $\|z\|_{S^m_1[M_m]} < 1$. Then by homogeneity we can write $z =
(a\otimes I_{M_m}) \cdot y \cdot (b\otimes I_{M_m})$ with $\|a\|_{S^m_2}
<1$, $\|y\|_{M_m(M_m)} <1$, $\|b\|_{S^m_2}<1$. Then $z=(z_{ij})$ with
$z_{ij}\in M_m$ defined by
$$z_{ij} = \sum_{k\ell} a_{ik}y_{k\ell}b_{\ell j}.$$
Hence we have
$$\leqalignno{|tr(z)| &= \left|\sum_{ij} \langle
e_i, z_{ij}e_j\rangle\right| &(1.10)\cr
&= \left|\sum_{k\ell} \langle\alpha_k, y_{k\ell} \beta_\ell\rangle
\right|}$$
where $\alpha_k, \beta_\ell \in \ell^m_2$ are defined by
$$\beta_\ell  =\sum_j b_{\ell j}e_j\quad \hbox{and}\quad \alpha_k = \sum_i
a_{ik}e_i.$$
Now (1.10) implies
$$\leqalignno{|tr(z)| &\le \|y\|_{M_m(M_m)} \left(\sum_k \|\alpha_k\|^2_2
\right)^{1/2} \left(\sum_{\ell}\|\beta_\ell\|^2_2\right)^{1/2}\cr
&\le \|y\|_{M_m(M_m)} \|a\|_{S^m_2} \|b\|_{S^m_2} < 1.}$$
By homogeneity, this concludes the proof of (1.9).

\vfill\eject
 
\n {\bf \S 2. Completely positive maps and regular maps.}

Let $K,H$ be Hilbert spaces. Consider (closed) subspaces $S\subset S_p(K)$
and $T\subset S_p(H)$, with $1\le p\le \infty$.  Let $E$ be an operator space. We will denote simply
by $S[M_n]$ (resp. $S[E]$) the space $S\otimes M_n$ equipped with the norm
induced by $S_p[K; M_n]$ (resp. the closure of $S\otimes E$ in $S_p[K;
E]$). We equip $S[M_n]$ and $S[E]$ with the operator space structures
induced by the {\it natural\/} ones on the spaces $S_p[K; M_n]$ and $S_p[K;
E]$. We will again refer to these o.s.s. on $S[M_n]$ and $S[E]$ as the {\it
natural\/} ones.

\proclaim Definition 2.1. In the preceding situation, we will say that a
linear map $u\colon \ S\to T$ is regular if $\sup\limits_{n\ge 1} \|u\otimes
I_{M_n}\|_{S[M_n]\to T[M_n]} <\infty$, and we denote
$$\|u\|_r = \sup_{n\ge 1}\|u\otimes I_{M_n}\|_{S[M_n]\to T[M_n]}.$$
We will denote by $B_r(S,T)$ the space of all regular maps from $S$ into
$T$. Equipped with the norm $\|\ \|_r$, this clearly is a Banach space.
(Actually, it can be equipped with the operator space structure of the
direct sum $\bigoplus\limits_{n\ge 1} cb(S[M_n], T[M_n])$.)

\n {\bf Remark.} In the case $p=\infty$  (resp. $p=1$), the spaces
${S[M_n]}$ and $ T[M_n]$ (resp.
${S_1[K,M_n]}$ and ${S_1[H,M_n]}$)
are tensor products associated 
with the minimal tensor product (resp. with 
the operator space version
of the projective tensor product). Therefore, in these two
 extreme cases
every $cb$ map $u$ is regular 
and satisfies $\|u\|_{cb}=\|u\|_r$.
Hence, in the case  $p=\infty$  
 we have  an isometric identity
$B_r(S,T)=cb(S,T)$, in particular
$$B_r(S_\infty(K),S_\infty(H))=cb(S_\infty(K),S_\infty(H)),$$
 and in the case  $p=1$, we have isometrically
$$B_r(S_1(K),S_1(H))=cb(S_1(K),S_1(H)).$$

Let  $u\colon \ S\to T$ be a regular linear map.  
Then for any operator space $E$, $u\otimes I_E$ extends to a bounded map
from $S[E]$ into $T[E]$ with
$$\|u\otimes I_E\|_{S[E]\to T[E]} \le \|u\|_r.\leqno (2.1)$$
This is rather easy to check using the following identity.

\n For any $a$ in $S[E]$
$$\|a\|_{S[E]} = \sup\|(I_S\otimes w)(a)\|_{S[M_m]}\leqno (2.2)$$
where the supremum runs over all $m \ge 1$ and all maps $w\colon \ E\to
M_m$ with $\|w\|_{cb}\le 1$. To prove (2.2), we can assume that $E =B(H)$
for some $H$ and that $S=S_p(K)$. 
The proof can then be completed using
Lemma~1.12 and Corollary~1.8 in [P3]. We leave the details to the reader.

\proclaim Proposition 2.2. Let $S,T$ be as above. Every regular map
$u\colon \ S\to T$ is $c.b.$ and satisfies
$$\|u\|_{cb}\le \|u\|_r.\leqno (2.3)$$

\pf This follows from Remark~2.4 in [P3]. There it is proved that the
$c.b.$ norm of $u$ can be written equivalently as follows
$$\|u\|_{cb} = \sup_{n\ge  1} \|I_{S^n_p}\otimes u\|_{S^n_p[S]\to
S^n_p[T]}.\leqno (2.4)$$
Actually, (2.4) is valid for a map $u\colon \ S\to T$ between arbitrary
operator spaces. But now if $S\subset S_p(K)$ and $T\subset S_p(H)$ we have
by (1.6) isometric identities
$$S^n_p[S] = S[S^n_p]\quad \hbox{and}\quad S^n_p[T] = T[S^n_p].$$
Hence applying (2.1) with $E = S^n_p$, we obtain (2.3).\qed

\proclaim Lemma 2.3. Consider a linear map $u\colon \ S_p\to S_p$. If $u$
is completely positive and bounded, then $u$ is regular and $\|u\|_r \le
\|u\|$. Moreover, $u\colon \ S_p\to S_p$ is regular iff its adjoint\break
$u^*\colon \ S_{p'}\to S_{p'}$ is regular and $\|u\|_r = \|u^*\|_r$.

\pf By a simple density argument, it suffices to prove this for $u\colon\
S^N_p\to S^N_p$ with $N\ge 1$ arbitrary. Then, if $u$ is $c.p.$, by
Stinespring's theorem (cf. [Pa, p. 53]) there is a finite set $y_1,\ldots, y_m$
in $M_N$ such that $u$ is of the form
$$u(x) = \sum^m_1 y_ixy^*_i,\qquad \forall x\in S^N_p.\leqno (2.5)$$
Assume $\|u\|=1$. Consider $a = (a_{ij})_{ij\le { N}}$ with $a_{ij}\in
M_n$ and $\|a\|_{S^N_p[M_n]} <1$. Then, by $(1.5)'$, we can write
$a = \alpha\cdot b\cdot \beta$ with $\alpha,\beta \in B_{S^N_{2p}}$ and
with $b=(b_{ij})$ in the unit ball of $M_N(M_n)$. Then we have
$$(u\otimes I_{M_n})(a) = \alpha_1\cdot \left(\matrix{b\cr &b&&\bigcirc\cr
&&\ddots\cr &\bigcirc&&b\cr}\right)\cdot \beta^*_1$$
where $\alpha_1 = (y_1\alpha,\ldots, y_m\alpha)$ and $\beta_1 =
(y_1\beta^*,\ldots, y_m\beta^*)$. By (1.5) this implies
$$\|(u\otimes I_{M_n})(a)\|_{S^N_p[M_n]} \le \|\alpha_1\|_{2p}
\|\beta^*_1\|_{2p},$$
whence
$$\eqalign{\|(u\otimes I_{M_n})(a)\|_{S^N_p[M_n]} &\le
({tr}(\alpha_1\alpha^*_1)^p)^{1/p} ({tr}(\beta_1\beta^*_1)^p)^{1/2p}\cr
&\le \left\|\sum y_i\alpha\alpha^*y^*_i\right\|^{1/2}_p \left\|\sum
y_i\beta^*\beta y^*_i\right\|^{1/2}_p\cr
&\le \|u(\alpha\alpha^*)\|^{1/2}_p \|u(\beta^*\beta)\|^{1/2}_p\cr
&\le \|u\|\le 1.}$$
This prove the first assertion.

\n We turn to  the second part. Using the fact that $S_p[M^*_n]^* =
S_{p'}[M_n]$ completely isometrically (cf. [P3]) and using (2.1) with
$E=M^*_n$ we find
$$\|u^*\|_r \le \|u\|_r.$$
Since $u = (u^*)^*$, the converse is obvious. This completes the proof.
\qed\medskip

\n {\bf Remark.} Let $S$ be a subspace of $\ell_p$. Since $\ell_p$ can be
viewed as embedded in $S_p$ via the diagonal matrices, there are two
notions of regularity for an operator $u\colon \ S\to \ell_p$. But it is
easy to check that $u\colon\ S\to \ell_p$ is regular in the Banach lattice
sense (see the introduction) iff it is regular in our sense (viewing $S$
and $\ell_p$ as embedded into $S_p$).

It can also be easily checked that $u\colon \ S\to \ell_p$ is completely
positive iff it is positive in the usual Banach lattice sense. This remark
clearly extends to more general $\ell_p$-spaces and explains why we
allowed  ourselves to use the same word ``regular'' as in the Banach
lattice case (while ``completely regular''  would have been tempting!).
\vfill\eject

\n {\bf \S 3. Main result and applications.}

To prove our main result on regular operators, we will need to describe the
predual of the space of regular operators on $S_p$. It will be simpler to
work in the finite dimensional case, so let $n$ and $m$ be fixed integers.

Consider $a$ in $S^n_p \otimes S^m_{p'}$ and assume that there are
$\alpha,\beta$ in $S^m_{2p'}$ and $y$ in $S^n_p[M_m]$ such that
$$a = (I\otimes \alpha)\cdot y \cdot (I\otimes \beta).\leqno (3.1)$$
We define
$$\rho_p(a) = \inf\{\|\alpha\|_{S^m_{2p'}} \|y\|_{S^n_p[M_m]}
\|\beta\|_{S^m_{2p'}}\},\leqno (3.2)$$
where the infimum runs over all possible representations of the form (3.1).

\n A more symmetric equivalent definition
 is as follows. Consider all
representations of $a$ of the form
$$a = (\gamma\otimes \alpha)\cdot g \cdot (\delta\otimes \beta),\leqno(3.1)'$$
with $\alpha,\ \beta$ as above, with $g$ in $M_n(M_m)$ and with
$\gamma$   and $\delta$ in $S^n_{2p}$. Then by $(1.5)'$ we 
clearly have
$$\rho_p(a) = \inf\{\|\gamma\|_{S^n_{2p}} \|\alpha\|_{S^m_{2p'}} \|g\|_{M_n(M_m)}
\|\beta\|_{S^m_{2p'}} \|\delta\|_{S^n_{2p}}\}.\leqno (3.2)'$$

We will prove below that $\rho_p$ is a norm. Surprisingly this is not so
obvious. The argument will be based on interpolation (as in [P2] for the
commutative case). Let $X^{nm}_p$ be the space $S^n_p\otimes S^m_{p'}$
equipped with the norm $\rho_p$. We will prove below that
$$(X^{nm}_p)^* = B_r(S^n_p, S^m_p)$$
isometrically and also that if $\theta =1/p$ we have isometrically
$$X^{nm}_p = (S^m_1[M_n], S^n_1[M_m])_\theta.\leqno (3.3)$$
Note that in (3.1) we have denoted by $I$ the identity matrix in $M_n$ and
the product appearing in (3.1) is the product in $M_n\otimes M_m$. However,
we will frequently need in the sequel to identify $M_n\otimes M_m$ and
$M_m\otimes M_n$ in the usual manner $(x\otimes y\to y\otimes x)$ so that
$y$ in (3.1) can be viewed alternatively as an element of $M_m\otimes M_n$
so that (3.1) becomes
$$a = (\alpha \otimes I) \cdot y \cdot (\beta\otimes I).\leqno (3.1)'$$
We warn the reader that we will use {\it both\/} ways to write (3.1) in the
sequel. The context will always make clear whether we work in $M_m \otimes
M_n$ or $M_n\otimes M_m$.  This identification is also used to give  a
meaning to the various interpolation theorems we consider in this section.
For instance it is used to view the couple appearing in (3.3) as a
compatible interpolation couple.

The proof is based on the matricial version of Szeg\"o's classical theorem
due to Masani-Wiener-Helson-Lowdenslager (see [He]) which can be stated as
follows (this form suffices for our purposes).

\n Let $w\colon \ {\bf T}\to M_m$ be a measurable matrix valued Lebesgue
integrable function. Assume that for some $\epsilon>0$ we have
$$\forall t\in {\bf T}\qquad w(t) \ge \epsilon I.$$
Then there is an analytic function $F\colon\ D\to M_m$ with entries in
$H^2$ such that the boundary values (= radial limits) satisfy
$$ \qquad F^*(z) F(z) = w(z)\quad {\rm a.e.\  on}\  \partial\Omega$$
and moreover such that $F(z)$ is invertible for all $z$ in $D$ and $z\to
F(z)^{-1}$ is bounded and analytic on $D$.

\proclaim Theorem 3.1. Consider $a$ in $S^n_p\otimes S^m_{p'}$. Let $\theta
= 1/p$. The following are equivalent.\medskip
\item{(i)} $\rho_p(a) <1$,
\item{(ii)} $a$ is in the open unit ball of the space $(S^m_1[M_n],
S^n_1[M_m])_\theta$, \medskip
\n Therefore $\rho_p$ is a norm and we have the isometric identity (3.3).

\pf Assume first $\rho_p(a)<1$. Then $a$ admits a representation (3.1)$'$ with
$\alpha,\beta$ (resp.\ $\gamma,\delta$) in the open unit ball of
$S^m_{2p'}$ (resp.\ $S^n_{2p}$) and $g$ in the open unit
ball of $M_n(M_m)$. Since
$$S^m_{2p'} = (S^m_2, S^m_\infty)_\theta \quad \hbox{and}\quad S^n_{2p} =
(S^n_\infty, S^n_2)_\theta,\leqno (3.4)$$
we may apply the classical multilinear interpolation theorem (cf.\ [BL,
p.~96] to the multilinear map
$$(\alpha, \beta, \gamma, \delta) \longrightarrow (\gamma\otimes \alpha)
\cdot g \cdot (\delta\otimes \beta).$$
We claim that this map is a contraction from
$$S^m_{2p'} \times S^m_{2p'}\times S^n_{2p}\times
S^n_{2p}\quad \hbox{into}\quad (S^m_1[M_n],
S^n_1[M_m])_\theta.$$ Indeed this is true if $\theta = 0$
and if $\theta=1$ (with $\theta=1/p$), hence by
interpolation this also holds for any intermediate
$0<\theta<1$. Hence we obtain that  (i) implies (ii).  

\n Conversely assume (ii). By definition of the complex interpolation space,
we can write $a=f(\theta)$ where $f$
 is an analytic function on the strip
$\Omega = \{0< {Re} z<1\}$,
 with values in $M_m \otimes M_n$, bounded and
 continuous on $\overline\Omega = \{0 \le {Re}
z\le 1\}$  and such that
$$\sup_{t\in {\reel}} \|f(it)\|_{S^m_1[M_n]} <1,\quad \sup_{t\in {\reel}}
\|f(1+it)\|_{S^n_1[M_m]} <1.$$
Let $\partial_0 = \{z\in {\comp}\mid {Re} z=0\}$ and $\partial_1 = \{z\in
{\comp}\mid {Re} z=1\}$. By a classical continuous selection argument, this
implies that we can find  bounded continuous functions $\alpha\colon \ \partial
\Omega\to M_m$ and $\beta\colon \ \partial\Omega\to M_m$ and $y\colon \
\partial\Omega\to M_n\otimes M_m$ such that
$$\forall  \ z\in \partial\Omega\qquad f(z) = (I\otimes \alpha(z))\cdot
y(z) \cdot (I\otimes \beta(z)),\leqno (3.5)$$
$$\sup_{z\in \partial_0}\|\alpha(z)\|_{S^m_2}<1,\quad \sup_{z\in
\partial_0}\|\beta(z)\|_{S^m_2}<1 \quad \hbox{and}\quad \sup_{z\in
\partial_0}\|y(z)\|_{M_m(M_n)}<1.\leqno (3.6)^1$$
$$\sup_{z\in \partial_1}\|\alpha(z)\|_{S^m_\infty} < 1, \quad \sup_{z\in
\partial_1}\|\beta(z)\|_{S^m_\infty} < 1\quad\hbox{and}\quad \sup_{z\in
\partial_1}\|y(z)\|_{S^n_1[M_m]} <1.\leqno (3.6)^2$$
Furthermore, we can write
$$\forall z\in \partial\Omega  \qquad y(z) = (\gamma(z)\otimes I)
\cdot  g(z)\cdot  (\delta(z)\otimes I) \leqno (3.5)'$$
with $\gamma(z), \delta(z)\in M_n$ and $g(z) \in M_n(M_m)$ such that
$$\leqalignno{&\sup_{z\in\partial_0} \|\gamma(z)\|_{S^n_\infty}<1,\quad
\sup_{z\in
\partial_0}\|\delta(z)\|_{S^n_\infty}<1, \hbox{ and } \sup_{z\in
\partial_0} \|g(z)\|_{M_n(M_m)} <1.& (3.6)^3\cr
&\sup_{z\in\partial_1}\|\gamma(z)\|_{S^n_2}<1,\quad \sup_{z\in
\partial_1}\|\delta(z)\|_{S^n_2}<1, \hbox{ and } \sup_{z\in
\partial_1}\|g(z)\|_{M_n(M_m)}<1.&(3.6)^4}$$
Note that for $z\in \partial_1$ (resp. $z\in\partial_0$) we can take
$\alpha(z)$ and $\beta(z)$ (resp. $\gamma(z)$ and $\delta(z)$) equal to a
multiple of the identity matrix. This yields
$$\forall z\in \partial\Omega \qquad f(z) = (\gamma(z)\otimes \alpha(z)) \cdot g(z)\cdot (\delta(z) \otimes
\beta(z)).\leqno(3.5)''$$
Of course the functions $\alpha, \beta, \gamma, \delta$ and $g$ are a
priori no longer (boundary values of) analytic functions, but we will
correct this using the matricial Szeg\"o theorem. We will now choose
$\varepsilon>0$ (small enough to be specified later) and we let for all $z$
in $\partial\Omega$
$$w_1(z) = \varepsilon I + \alpha(z) \alpha(z)^*\qquad v_1(z) = \varepsilon
I + \gamma(z) \gamma(z)^*$$
and
$$w_2(z) = \varepsilon I + \beta(z)^* \beta(z)\qquad v_2(z) = \varepsilon I
+ \delta(z)^* \delta(z).$$
Then by the matricial version of Szeg\"o's theorem (and by a conformal
mapping argument), there are $M_m$-valued (resp.\ $M_n$-valued) analytic
functions $z\to \tilde\alpha(z)$ and $z\to \tilde\beta(z)$ (resp.\ $z\to
\tilde\gamma(z)$ and $z\to \tilde\delta(z)$) such that $z\to
\tilde\alpha(z)^{-1}$ and $z\to \tilde\beta(z)^{-1}$ (resp.\ $z\to
\tilde\gamma(z)^{-1}$ and $z\to \tilde\delta(z)^{-1}$) are well defined,
bounded and analytic in $\Omega$ and such that for almost all $z$ in
$\partial\Omega$
$$w_1(z) =\tilde\alpha(z) \tilde\alpha(z)^*,\ v_1(z) =\tilde\gamma(z)
\tilde\gamma(z)^*, \ w_2(z) = \tilde\beta(z)^*
\tilde\beta(z), \ v_2(z) =
\tilde\delta(z)^*\tilde\delta(z).\leqno (3.7)$$ 
We can
then write for all $z$ in $\Omega$ $$f(z) =
(\tilde\gamma(z)\otimes\tilde \alpha(z))\cdot \tilde g(z)
\cdot (\tilde \delta(z) \otimes\tilde \beta(z))$$ where 
 we have set
$$\tilde g(z)  =
(\tilde\gamma(z)^{-1}\otimes \tilde \alpha(z)^{-1}) \cdot
f(z) \cdot (\tilde\delta(z)^{-1} \otimes \tilde
\beta(z)^{-1}) . \leqno(3.8)$$
Note that by (3.8) $\tilde g$ is bounded and analytic  in $\Omega$ and its
boundary values satisfy (by
 (3.5)$''$) a.e. on $\partial\Omega$
$$\tilde g(z)= u(z)\cdot g(z) \cdot v(z),$$
where $u(z) = \tilde \gamma(z)^{-1} \gamma(z) \otimes \tilde\alpha(z)^{-1}
\alpha(z)$ and $v(z) = \delta(z)\tilde \delta(z)^{-1}\otimes \beta(z)
\tilde\beta(z)^{-1}$. But by (3.7) for almost all $z$ in
$\partial\Omega$ we have $$\tilde\alpha(z)
\tilde\alpha(z)^* \ge \alpha(z)\alpha(z)^* \hbox{ and }
\tilde\gamma(z)\tilde\gamma(z)^* \ge \gamma(z) \gamma(z)^*$$
hence $\|u(z)\|_{M_n\otimes_{\rm min} M_m} = \|u(z)\|_{M_n(M_m)}\le 1$.
Similarly, $\|v(z)\|_{M_n(M_m)} \le 1$. Therefore,  the
boundary values of $\tilde g$ satisfy the same bounds (3.6)$^3$ and $(3.6)^4$ as $g$ on
$\partial\Omega = \partial_0 \cup\partial_1$. Since $\tilde g$ is
bounded and analytic,
this implies by the maximum principle
$$\forall \ z\in \Omega\qquad \qquad \|\tilde g(z)\|_{M_n(M_m)}\le 1.$$
On the other hand, if we choose $\varepsilon$ small enough we can by (3.7)
guarantee
that $\tilde\alpha, \tilde \beta, \tilde\gamma, \tilde\delta$ still satisfy
the estimates (3.6)$^1$, (3.6)$^2$, (3.6)$^3$ and (3.6)$^4$. We then
obtain
$$\|\tilde\alpha(\theta)\|_{(S^m_2,S^m_\infty)_\theta}<1,\ 
\|\tilde\gamma(\theta\|_{(S^n_\infty, S^n_2)_\theta}<1\
\hbox{ and }\
\|\tilde\beta(\theta)\|_{(S^m_2,S^m_\infty)_\theta} <1,\ 
\|\tilde \delta(\theta)\|_{(S^n_\infty,
S^n_2)_\theta}<1.$$ Then using (3.4) and recalling
(3.2)$'$ we conclude that $\rho_p(a)<1$ follows from the identity
$$a=f(\theta) =
(\tilde\gamma(\theta)\otimes \tilde\alpha(\theta)) \cdot
\tilde g(\theta) \cdot (\tilde\delta(\theta)\otimes
\tilde\beta(\theta)).$$  This
concludes the proof (ii) $\Rightarrow$ (i).
\qed

\proclaim Theorem 3.2. In the situation of Theorem~3.1, we have
isometrically 
$$(X^{nm}_p)^* = B_r(S^n_p, S^m_p).$$

\pf By Theorem~3.1 we have
$$(X^{nm}_p)^* = (S^m_1[M_n]^*, \quad S^n_1[M_m]^*)_\theta.$$

Observe the isometric identities
$$\eqalign{S^m_1[M_n]^* &\simeq cb(M_n,M_m) = B_r(M_n,M_m)\cr
S^n_1[M_m]^* &\simeq cb(S^n_1,S^m_1) = B_r(S^n_1,S^m_1).}$$
Using (1.4) with $E=M_k$ with $k$ arbitrary, it is easy to show that (if
$p=1/\theta$) we have a norm one inclusion
$$(B_r(M_n,M_m), \quad B_r(S^n_1, S^m_1))_\theta \subset B_r(S^n_p,
S^m_p).$$
This shows that we have a norm one inclusion
$$(X^{nm}_p)^* \subset B_r(S^n_p, S^m_p).$$
To prove the converse it suffices to prove the following claim:\ for any
$u$ in $B_r(S^n_p, S^m_p)$ and any $a$ in $X^{nm}_p$ we have
$$|\langle u,a\rangle| \le \rho_p(a)\|u\|_r.$$
Let us verify this. By homogeneity we may assume $\rho_p(a)<1$. Then we can
assume that (3.1) holds with
$$\|\alpha\|_{S^m_{2p'}}<1,\quad \|y\|_{S^n_p[M_m]} < 1, \quad
\|\beta\|_{S^m_{2p'}} < 1.$$
Let $z = (u\otimes I_{M_m})(y)$. Note $z\in S^m_p[M_m]$. We have by
definition of $\|u\|_r$
$$\|z\|_{S^m_p[M_m]} \le\|u\|_r \|y\|_{S^n_p[M_m]} < \|u\|_r.$$
Let $\widehat z = (^t\alpha\otimes I_{M_m}) \cdot z \cdot (^t\beta\otimes
I_{M_m})$. We have by $(1.5)'$
$$\|\widehat z\|_{S^m_1[M_m]} \le \|\alpha\|_{2p'} \|z\|_{S^m_p[M_m]}
\|\beta\|_{2p'} <\|u\|_r.$$
On the other hand by (1.8)
$$\langle u,a\rangle = tr((I_{M_m}\otimes \alpha)z (I_{M_m}\otimes \beta))
= tr(\widehat z)$$
hence by (1.9)
$$|\langle u,a\rangle| \le \|\widehat z\|_{S^m_1[M_m]} < \|u\|_r.$$
This concludes the proof of our claim.\qed

\proclaim Corollary 3.3. Let $n,m$ be arbitrary integers, $0\le \theta \le
1$, $p=1/\theta$. Then we have isometric identities
$$\leqalignno{&(cb(S^n_\infty, S^m_\infty), \quad cb(S^n_1,S^m_1))_\theta =
B_r(S^n_p, S^m_p),&(3.9)'\cr
&(cb(S_\infty, B(\ell_2)), \quad cb(S_1,S_1))^\theta = B_r(S_p,
S_p).&(3.9)''}$$

\pf Let $X^{nm}_\infty = S^m_1[M_n]$ and
$X^{nm}_1 = S^n_1[M_m]$. Then by Theorem~3.1 we have isometrically
$$X^{nm}_p = (X^{nm}_\infty, X^{nm}_1)_\theta.\leqno (3.10)$$
On the other hand, by our definitions and known estimates (cf. [BP, ER2,
ER6])
$$\eqalignno{X^{nm*}_1 &= (S^n_1)^* \otimes_{\rm min}(M_m)^* \simeq
cb(S^n_1, S^m_1)\cr
\noalign{\hbox{and}}
X^{nm*}_\infty &= (S^m_1)^* \otimes_{\rm min} (M_n)^* \simeq cb(M_n,
M_m).}$$
Hence $(3.9)'$ follows from (3.10) by duality, using Theorem 3.2. By an entirely elementary
approximation argument (left to the reader) we can obtain
$(3.9)''$.\qed\medskip

\n {\bf Proof of Theorem 0.1.} We simply combine $(3.9)''$ with our earlier
result Theorem~0.2. \qed
We may clearly extend the definition (3.2) to the infinite dimensional
case. For any $a$ in $S_p(H)\otimes S_{p'}(K)$ we consider all possible
factorizations of $a$ of the form
$$a = (I\otimes \alpha) \cdot y \cdot (I\otimes \beta)$$
with $\alpha,\beta\in S_{2p'}(K)\quad \hbox{and}\quad y \in S_p[H;
S_\infty(K)]$. We then let
$$\rho_p(a) = \inf\{\|\alpha\|_{2p'} \|y\|_{S_p[H; S_\infty(K)]}
\|\beta\|_{2p'}\}$$
and we define $X^{H,K}_p$ as the completion of $S_p(H) \otimes S_{p'}(K)$
for this norm.

Actually we will work only with $H=K=\ell_2$ and in that case  we set $X_p
= X^{H,K}_p$. Alternatively, we could use the norms defined in (3.2), form
the inductive limit $\bigcup\limits_{n,m} S^n_p \otimes S^m_{p'}$ and
define $\rho_p$ as the resulting completion. In any case, we have clearly
by (3.3) and Theorem~3.1 the following result

\proclaim Corollary 3.4. The following are isometric identities
$(\theta=1/p,\  1<\theta<1$)
$$X^{H,K}_p = (S_1[K; S_\infty(H)], \quad S_1[H; S_\infty(K)])_\theta$$
and
$$(X^{H,K}_p)^* = B_r(S_p(H), S_p(K)).$$

We will now prove an extension theorem,
 which in the commutative case with
 $p=1$ goes back to M.~L\'evy [L\'e]. See [P2] for the
commutative case with $1\le p<\infty$ arbitrary.

\proclaim Theorem 3.5. Let $1\le p <\infty$. Let $S\subset S_p$ be a closed
subspace. Then any regular operator $u\colon \ S\to S_p$ admits a regular
extension $\tilde u\colon \ S_p\to S_p$ with $\|\tilde u\|_r = \|u\|_r$.

\pf This follows from the Hahn-Banach theorem and a special property of the
norm of the predual of $B_r(S_p,S_p)$. It clearly suffices to consider the
case of a map $u\colon \ S\to S^m_p$. Assume $\|u\|_r\le 1$. We first claim
that it suffices to prove the following:\ for all $v$ in $S\otimes
S^m_{p'}$ the element $\tilde v \in S_p \otimes S^m_{p'}$ associated to $v$
by the inclusion $S\to S_p$ satisfies
$$|\langle u,v\rangle| \le \|\tilde v\|_{X_p}.\leqno (3.11)$$
Here for simplicity we again denote by $X_p$ the space $X^{HK}_p$ associated
to $H=\ell_2$, $K=\ell^m_2$. Indeed, if (3.11) holds the Hahn-Banach
theorem provides an element $\tilde u \in X^*_p$ with\break $\|\tilde
u\|_{X^*_p}\le 1$ such that $\langle \tilde u,v\rangle = \langle
u,v\rangle$ for all $v$ in $S\otimes S^m_{p'}$. Since this last condition
means that $\tilde u$ extends $u$ and since $X^*_p = B_r(S_p,S^m_p)$ we
obtain the desired conclusion.

\n Hence it suffices to show (3.11). Assume $\|\tilde v\|_{X_p} <1$. Then
there is a factorization $\tilde v = (I\otimes \alpha) y(I\otimes \beta)$
with $\|\alpha\|_{S^m_{2p'}} <1$, $\|\beta\|_{S^m_{2p'}} <1$ and
$\|y\|_{S_p[M_m]} < 1$. By a perturbation argument, we may clearly assume
that  $\alpha$ and $\beta$ are invertible. Then we have
$$y = (I\otimes \alpha^{-1}) \tilde v(I\otimes \beta^{-1})$$
which shows that $y \in S\otimes M_m$ since $\tilde v$ comes from $v\in
S\otimes S^m_{p'}$.
 Hence we can write $v = (I\otimes \alpha) y(I\otimes \beta)$.

\n The rest of the proof is a variant of the proof of
Theorem~3.2. Let $z = (u\otimes I_{M_m})(y)$, note $z\in
S^m_p[M_m]$, and let $\widehat z = (^t\alpha \otimes I)
z(^t\beta\otimes I)$. Then by $(1.5)'$ we have $$\|\widehat
z\|_{S^m_1[M_m]} < \|z\|_{S^m_p[M_m]} < \|u\|_r.\leqno
(3.12)$$ Hence we conclude using (1.8) and (1.9)
$$|\langle u,v\rangle| = |tr((I\otimes \alpha) z(I\otimes \beta))| =
|tr(\widehat z)| \le \|\widehat z\| _{S^m_1[M_m]}.$$
This shows that (3.12) implies (3.11).\qed

To illustrate the preceding results, we consider the subspace $T_p\subset
S_p$ of all the upper triangular matrices. This space is often regarded as
a noncommutative analogue of $H^p$.

Using the results of [P4], it is not hard to show that we have
isomorphically if $1<p<\infty$ and $\theta =1/p$
$$T_p[S_\infty] = (T_\infty[S_\infty], T_1[S_\infty])_\theta.$$
Hence by repeating mutatis mutandis the proof of [P2] in the case of $H^p$,
we find

\proclaim  Corollary 3.6. We have (if $1<p<\infty$ and $\theta = 1/p$) an
isomorphic identity
$$B_r(T_p,S_p) = (B_r(T_\infty,S_\infty), \quad B_r(T_1,S_1))^\theta.$$

Following the framework of [P3], it is easy to extend the preceding results
to the case when $S_p$ or $S^n_p$ is replaced by a non-commutative
$L_p$-space associated to
 a {\it hyperfinite\/} von~Neumann algebra $M$. Recall
that a von Neumann algebra $M$ is called 
hyperfinite if it is the $\sigma(M,M_*)$-closure of the
union of an increasing net of finite dimensional
subalgebras. The
extension of the definition of a regular operator
 is immediate:

\n Let $M$
(resp. $N$) be a hyperfinite von~Neumann algebra equipped with a faithful
normal semi-finite trace $\varphi$ (resp. $\psi$).
We will denote by $L_p(\varphi)$ (resp. $L_p(\psi)$) the associated
non-commutative $L_p$-space for $1\le p <\infty$. (Note that it is natural
to identify $L_\infty(\varphi)$ with $M$ and $L_\infty(\psi)$ with $N$.)
These spaces are equipped with their {\it natural o.s.s.\/} as explained in
section~1. Similarly, if $E$ in an operator space, the space
$L_p(\varphi;E)$ is defined by interpolation as in (1.4) and it is equipped
with the natural o.s.s. defined by (1.4).

\n Consider then a closed subspace
$$S\subset L_p(\varphi) \quad (\hbox{resp. } T\subset L_p(\psi)).$$
We will denote by $S[E]$ the closure of $S\otimes E$ in the space
$L_p(\varphi; E)$. We can then define a regular map $u\colon \ S\to T$
exactly as in definition~2.1.

\n We again denote by $B_r(S,T)$ the space of all regular maps
from $S$ into $T$ equipped with the norm $\|\ \|_r$.
It is then very easy to adapt the
preceding proofs to obtain the following two statements. We
leave the details to the reader.

\proclaim Theorem 3.7. Let $M$ (resp. $N$) be a hyperfinite von~Neumann algebra
equipped with a semi-finite faithful normal trace $\varphi$
(resp. $\psi$). Let $1<p<\infty$ and
$\theta =1/p$. We have an isometric identity
$$(cb(M,N), cb(L_1(\varphi), L_1(\psi)))^\theta = B_r(L_p(\varphi),
L_p(\psi)).$$
Moreover the space $B_r(L_p(\varphi), L_p(\psi))$ coincides with the set
of all maps $u\colon \ L_p(\varphi)\to L_p(\psi)$ which are linear
combinations of completely positive, bounded maps from
$L_p(\varphi)$ to $L_p(\psi)$.

\proclaim Corollary 3.8. Let $1\le p<\infty$.
 Let $(M,\varphi)$ and $(N,\psi)$ be as in Theorem 3.7.
 Let $S\subset L_p(\varphi)$
be a closed subspace.
Then any regular operator $u\colon\  S\to L_p(\psi)$ admits an extension
$\widetilde u\colon \ L_p(\varphi)\to L_p(\psi)$ with $\|\widetilde u\|_r =
\|u\|_r$.

\n {\bf Remark.} For the case $p=1$ in corollary 3.8,
 it is useful
to remind the reader that (by e.g. [Ta, p. 126-127]) 
there is a
norm one, completely positive and completely contractive
projection from $N^*=L_1(\psi)^{**}$ onto $L_1(\psi)$.

\n {\bf Remark.} Note that in the case $p=\infty$, the preceding two statements
reduce to the well known decomposition and extension
 properties of a completely
bounded map with values into 
a hyperfinite (=injective, by Connes's well known
results [Co]) von Neumann algebra, which follow from Wittstock's theorem (cf. [Pa,
p.100-107]). Note that these properties are only true in the injective case
(cf. [H1]).

\vfill\eject

 \centerline{\bf References}

 \item{[BL]} J. Bergh and J. L\"ofstr\"om. Interpolation spaces. An
 introduction. Springer Verlag, New York. 1976.

 \item{[B1]} D. Blecher. 
Tensor products of operator spaces II. 
 Canadian J. Math. 44 (1992) 75-90.
 
 \item{[B2]}  $\underline{\hskip1.5in}$. The standard
 dual of an operator space.
  Pacific J. Math. 153 (1992) 15-30.

\item{[B3]}  $\underline{\hskip1.5in}$. Generalizing Grothendieck's
program. in"Function spaces", Edited by K.Jarosz, Lecture
Notes in Pure and Applied Math. vol.136, Marcel Dekker,
1992.
 
 \item{[BP]} D. Blecher and V. Paulsen. Tensor products of operator spaces.
 J. Funct. Anal. 99 (1991) 262-292.
 
 \item{[BS]} D. Blecher and R. Smith. The dual of the Haagerup tensor
 product. Journal London Math. Soc. 45 (1992) 126-144.

 \item{[Ca]} A. Calder\'on. Intermediate spaces and interpolation, the
 complex method. Studia Math. 24 (1964) 113-190.

\item{[Co]} A. Connes. Classification of injective
factors, Cases
$II_1,II_\infty,III_\lambda,\lambda\neq 1$. Ann.
Math. 104 (1976) 73-116.

 \item{[ER1]} E. Effros and Z.J. Ruan. On matricially normed spaces. Pacific
 J. Math. 132 (1988) 243-264.
 
 \item{[ER2]} $\underline{\hskip1.5in}$. A new approach to operators spaces.
 Canadian Math. Bull.
34 (1991) 329-337.

 \item{[ER3]} $\underline{\hskip1.5in}$. On the abstract 
characterization of
 operator spaces. Proc. A.M.S. To appear.
 
 \item{[ER4]} $\underline{\hskip1.5in}$. Self duality for the Haagerup
 tensor product and Hilbert space factorization.  J. Funct. Anal. 100
(1991) 257-284.
 
 \item{[ER5]} $\underline{\hskip1.5in}$. Recent development in operator
 spaces.

 \item{[ER6]} $\underline{\hskip1.5in}$. Mapping spaces
and liftings for operator spaces. (Preprint) Proc. London
Math. Soc. To appear.

\item{[ER7]}  $\underline{\hskip1.5in}$. The
Grothendieck-Pietsch and Dvoretzky-Rogers Theorems for
operator spaces. (Preprint 1991) J. Funct. Anal. To
appear.

\item{[ER8]} $\underline{\hskip1.5in}$. On approximation
properties for operator spaces, International J. Math. 1
(1990) 163-187.

\item{[G]} A. Grothendieck.  R\'esum\'e de la th\'eorie
m\'etrique des produits tensoriels topologiques.  Boll..
Soc. Mat. S$\tilde{a}$o-Paulo  8 (1956), 1-79.

 \item{[H1]} U. Haagerup. Injectivity and decomposition of completely
 bounded maps in ``Operator algebras and their connection with Topology and
 Ergodic Theory''. Springer Lecture Notes in Math. 1132
(1985) 170-222.
 
 \item{[H2]} $\underline{\hskip1.5in}$. $L^p$-spaces
associated with an arbitrary von Neumann algebra.
Alg\`ebres d'op\'erateurs et leurs applications en
physique math\'ematique. (Colloque CNRS, Marseille, juin
1977) Editions du CNRS, Paris 1979.

 \item{[H3]} $\underline{\hskip1.5in}$. Decomposition of completely bounded
 maps on operator algebras. Unpublished manuscript. Sept. 1980.
 
\item{[HeP]} A. Hess and G. Pisier, On the
$K_t$-functional for the couple $B(L_1,L_1), B(L_\infty,
L_\infty)$). Quarterly J. Math. Oxford (submitted).

\item {[KR]} R. Kadison  and J. Ringrose.    Fundamentals
of the theory of operator algebras, Vol. II, Advanced
Theory,    Academic Press, New-York      1986.

\item{[Ko]} H. Kosaki. Applications of the complex
interpolation method to a von Neumann algebra: non
commutative $L^p$-spaces. J. Funct. Anal. 56 (1984) 29-78.

\item{[N]} E. Nelson. Notes on non-commutative
integration. J. Funct. Anal. 15 (1974) 103-116. 
  
 \item {[L\'e]} M. L\'evy. Prolongement d'un op\'erateur
d'un sous-espace de $L^1(\mu)$ dans $L^1(\nu)$. S\'eminaire
d'Analyse Fonctionnelle 1979-1980. Expos\'e 5. Ecole
Polytechnique.Palaiseau.

 \item{[Pa]} V. Paulsen.  Completely
bounded maps and dilations. Pitman Research Notes 146.
Pitman Longman (Wiley) 1986.

 \item{[P1]} G. Pisier. The operator Hilbert space $OH$,
complex interpolation and tensor norms. To appear.

\item {[P2]} $\underline{\hskip1.5in}$. Complex
interpolation and regular operators between Banach
lattices. Arch. der Mat. (Basel) To appear.

\item{[P3]} $\underline{\hskip1.5in}$.{  Noncommutative
vector valued $L_p$-spaces and completely
$p$-summing maps.} Preprint. To appear.

\item{[P4]} $\underline{\hskip1.5in}$. Interpolation of
$H^p$-spaces and noncommutative generalizations I. Pacific
J. Math. 155 (1992) 341-368.

\item{[P5]} $\underline{\hskip1.5in}$. Espaces $L_p$ non
commutatifs \`a valeurs vectorielles et applications
$p$-sommantes. C. R. Acad. Sci. Paris, 316 (1993)
1055-1060.

\item{[Ru]} Z.J. Ruan. Subspaces of $C^*$-algebras. J. Funct. Anal. 76
 (1988) 217-230.

\item{[Te1]} M. Terp. Interpolation spaces between a von
Neumann algebra and its predual. J. Operator Th.
8 (1982) 327-360.

\item{[Te2]} $\underline{\hskip1.5in}$. $L^p$-spaces associated with von
Neumann algebras. Preprint. Copenhagen University. June
1981.

\vskip12pt

\vskip12pt

Texas A. and M. University

College Station, TX 77843, U. S. A.

and

Universit\'e Paris 6

Equipe d'Analyse, Bo\^\i te 186,
 
75252 Paris Cedex 05, France

 \end